\documentclass[10pt,a4paper]{article}
\makeatletter
\usepackage[dvipsnames]{xcolor}
\usepackage{lipsum}
\usepackage{tipa}
\usepackage{tikz}
\usepackage{tikz-cd}
\usetikzlibrary{matrix,arrows,decorations.pathmorphing}
\newcounter{magicrownumbers}

\usepackage[left=3cm,top=3cm,bottom=3cm,right=3cm,head=1cm,foot=1cm]{geometry}
\usepackage[pdftex,bookmarks=true,bookmarksnumbered=true]{hyperref}
\hypersetup{linkbordercolor=ForestGreen!25, citebordercolor=Goldenrod!60}
\DeclareMathAlphabet{\mathpzc}{OT1}{pzc}{m}{it}
\usepackage{graphicx}
\usepackage{stmaryrd}
\usepackage{twoopt}
\usepackage{amssymb}
\usepackage{amsmath}
\usepackage{amsthm}
\usepackage{mathabx} 
\usepackage{amsfonts}
\usepackage{amscd}
\usepackage{url}
\usepackage[all]{xy}

\newcommand{\hil}{\mathcal{H}}

\newcommand{\5}{\hspace{0,5cm}}
\newcommand{\3}{\vspace{0,3cm}}
\newcommand{\hor}{\vspace{0,1cm}}
\newcommand{\la}{\langle}
\newcommand{\ra}{\rangle}

\newcommand{\fp}[2]{\mspace{1mu} {}_{#1} \mspace{-5mu} \times_{#2}}

\newcommand{\gpdhom}[7]{\begin{tikzpicture}
\path (3,1) node (a) {${#1}$} (0,0) node (b) {${#2}$} (2,0) node (c) {${#3}$} (4,0) node (d) {${#4}$} (6,0) node (e) {${#5}$}; 
\begin{scope}
\draw[->] (a) -- node [left] {${#6}$} (c) ; 
\draw[->] (a) -- node [right] {${#7}$} (d); 
\draw[->] (0.8,0.05) -- (1.5,0.05) ;
\draw[->] (0.8,-0.05) -- (1.5,-0.05) ;
\draw[->] (5.2,0.05) -- (4.5,0.05) ;
\draw[->] (5.2,-0.05) -- (4.5,-0.05) ;
\draw[->] (0,0.4) arc (180:90:64pt and 18pt);
\draw[->] (6,0.4) arc (0:90:64pt and 18pt);
\end{scope}
\end{tikzpicture}}

\title{Morita Equivalence and Spectral Triples on Noncommutative Orbifolds}
\date{}
\author{Antti J. Harju\footnote{harjuaj@gmail.com, Helsinki U / QMU London}}
\begin{document}
\maketitle

\begin{abstract}
Let $G$ be a finite group. Noncommutative geometry of unital $G$-algebras is studied. A geometric structure is determined by a spectral triple on the crossed product algebra associated with the group action. This structure is to be viewed as a representative of a noncommutative orbifold. Based on a study of classical orbifold groupoids, a Morita equivalence for the crossed product spectral triples is developed. Noncommutative orbifolds are Morita equivalence classes of the crossed product spectral triples. As a special case of this Morita theory one can study freeness of the $G$-action on the noncommutative level. In the case of a free action, the crossed product formalism reduced to the usual spectral triple formalism on the algebra of $G$-invariant functions. \hor


\end{abstract}

\section*{Introduction}

We take up the task to develop a noncommutative geometric model for unital algebras subject to a finite group action. The differential geometric objects relevant to this theory are global action orbifolds. These objects may be naturally represented as Lie groupoids, or more precisely, proper {\' e}tale Lie groupoids, \cite{Moe02}. It was shown in the reference \cite{Har14a} that with any effective compact proper {\' e}tale groupoid with a spin structure, one can associate a Dirac spectral triple on the smooth groupoid convolution algebra. So, from this viewpoint, the generalization of noncommutative geometry to $G$-algebras is straightforward. A mayor difference between the manifold theory and the orbifold theory is the implementation of an equivalence relation. Namely, orbifolds (considered as Lie groupoids) should be equipped with the geometric Morita equivalence relation which is essentially weaker relation than the $G$-equivariant diffeomorphism. The geometric Morita equivalence leaves the orbit space of the groupoid invariant, as well as the isotropy types associated with the action. So, one should view an orbifold as a Morita equivalence class of its representative groupoid since the different Morita classes are merely different geometric models for the same orbit space. A particularly important case is an orbifold that is modelled as a manifold subject to a free action of a finite group. This orbifold is Morita equivalent to the unit groupoid of its orbit space. So, the orbifolds subject to free actions should be viewed as manifold objects in the category of Lie groupoids. This is how the geometric Morita equivalent provides a tool to separate smooth objects from those with singularities. 

Suppose that a finite group $G$ acts on the unital complex algebra $A$. This results in two unital associative algebras of interest: the crossed product algebra $G \ltimes A$ and the invariant subalgebra $A^G$. The philosophy of this project is to view the algebra $G \ltimes A$ as an equivariant space whereas a spectral triple on $G \ltimes A$ is viewed as a representative of a noncommutative orbifold. We are also assuming that the $G$-action on the Hilbert space commutes with the Dirac operator which results in a spectral triple on the invariant subalgebra $A^G$. This level should be viewed as a noncommutative model for the orbit space under the group action. The spectral triple on $A^G$ captures the metric aspects of this theory which is demonstrated in the appendix of this document. These spectral triples also have an important role in the Morita theory introduced in Section 1. 

The main goal of this work is to implement a Morita equivalence for spectral triples on noncommutative orbifolds. We shall follow the standard philosophy of noncommutative geometry and translate the geometric Morita equivalence of orbifolds (viewed as Lie groupoids) into the operator theoretic language of spectral triples. Suppose that $G \ltimes X$ is a representative of a global action orbifold and that there is a Dirac spectral triple on its smooth crossed product algebra $(G \ltimes C^{\infty}(X), L^2(F_{\Sigma}) ,\eth)$. If $K \ltimes Y$ is any other global action orbifold that is Morita equivalent to  $G \ltimes X$, and if the Morita equivalence $\phi$ has been fixed, then there is the induced spectral triple $(K \ltimes C^{\infty}(Y), L^2(\phi_{\#} F_{\Sigma}),\phi_{\#} \eth)$ which was found in \cite{Har14b}. The induced spectral triple is defined in purely geometric means. It is known that once the geometric Morita equivalence $\phi$ has been fixed, then $\phi$ gives rise to a pre-$C^*$-algebra $K \ltimes C^{\infty}(Y)$ - $G \ltimes C^{\infty}(X)$ bimodule which describes the Morita equivalence on the algebraic level, \cite{MRW89}. This bimodule has a completion to a Morita bimodule for the $C^*$-algebra completions. In this manuscript we shall see that the assignment $L^2(F_{\Sigma}) \mapsto L^2(\phi_{\#}F_{\Sigma})$  is, up to a unitary equivalence, the map that induced an isomorphism of the pre-$C^*$-algebra representation categories, \cite{Rie74}. It will also be shown that the assignment $\eth \mapsto \phi_{\#} \eth$ will preserve the K-homology class of the underlying Fredholm module, i.e. in terms of operator KK-theory, the approximate sign of $\phi_{\#} \eth$ is a connection on $L^2(\phi_{\#}F_{\Sigma})$ for the approximate sign of $\eth$. However, these properties alone are too weak to provide axioms for a realistic model for a Morita equivalence of noncommutative orbifolds because the underlying Fredholm module loses information of the Dirac spectrum. The Dirac spectrum carries information of the metric associated with the riemannian structure. More precisely, in the appendix we shall see that in the case of global action orbifolds,  the Dirac spectrum captures the geodesic lengths already on the level of spectral triples on the algebras of smooth invariant functions $C^{\infty}(X)^G$. This is indeed fully compatible with the intuition: the geodesic length is a $G$-invariant property, and such properties are Morita invariant. Therefore, we need to introduce additional axiom for the Morita equivalence and require that a Morita equivalence operates as a unitary equivalence on the level of invariant spectral triples. This holds in the case of geometric orbifolds, \cite{Har14b}. 

\3 \noindent \textbf{Notation.} This work is a part of the project consisting of \cite{Har14a}, \cite{Har14b} and this manuscript. The previous parts are purely geometric studies of spectral triples on proper {\' e}tale groupoids. Although we do not need this generality in this manuscript, we shall continue to work with the Lie groupoid notation for the sake of coherency. If $X$ is an equivariant manifold subject to an action of a finite group $G$, then we consider this system as the action Lie groupoid $G \times X \rightrightarrows X$, which will be denoted by $G \ltimes X$. The source and target maps are given by $s(g, x) = x$ and $t(g, x) = g \cdot x$ for all $g \in G$ and $x \in X$. We also use $t^{-1}(x) = (G \ltimes X)^x$ and $s^{-1}(x) = (G \ltimes X)_x$

In the analysis, we shall identify the smooth crossed product algebras $G \ltimes C^{\infty}(X)$ with the groupoid convolution algebras $C^{\infty}(G \ltimes X)$. A smooth left Haar system $\mu = \{ \mu^x: x \in X\}$ in $G \ltimes X$ is a collection of measures in $G \times X$  parametrized on $X$, so that 
\begin{quote}
\textbf{1.} The support of $\mu^x$ is $(G \ltimes X)^x$. 

\textbf{2.} If $f \in C^{\infty}(G \ltimes X)$, then $x \mapsto \int f d \mu^x$ is smooth. 

\textbf{3.} The measures are left-invariant: $\int f(\sigma \tau) d \mu^{s(\sigma)}(\tau) = \int f(\tau) d \mu^{t(\sigma)}(\tau)$. 
\end{quote}
Here we are concerned with the action groupoids only. In this case we fix the standard Haar measure $\lambda$ in the finite group $G$ (the counting measure) and then define a Haar system in $G \ltimes X$ by setting
\begin{eqnarray*}
\mu = \{\mu^x = \lambda : x \in X\}.
\end{eqnarray*}
The groupoid convolution $C^{\infty}(G \ltimes X)$ is the unital algebra equipped with the product
\begin{eqnarray*}
(f \cdot g)(\sigma) &=& \int_{\tau \in (G \ltimes X)_{s(\sigma)}} f(\sigma \cdot \tau^{-1}) g(\tau) d \mu^{s(\sigma)}(\tau) \\
&=& \frac{1}{\# G} \sum_{\tau \in (G \ltimes X)_{s(\sigma)}} f(\sigma \cdot \tau^{-1}) g(\tau) 
\end{eqnarray*}
where $\# G$ denotes the  number of elements in $G$. The unit is the function $\textbf{1}_e$ which is the indicator function of the subset $\{e\} \times X$ in $G \times X$ where $e$ is the unit element of $G$. We shall drop the measure terms $d \mu^{s(\sigma)}$ in the notation. The algebra $C^{\infty}(G \ltimes X)$ can be equipped with the universal $C^*$-norm and the groupoid $C^*$-algebra $C^{\star}(G \ltimes X)$ is obtained by taking the completion in this norm \cite{Ren80}, \cite{Con94}.  

\3 \noindent \textbf{Convention.} An action Lie groupoid $G \ltimes X$ is called a compact action orbifold if $G$ is a finite group, $X$ is a compact connected  smooth manifold and the $G$ action on $X$ is effective. This is a nonstandard notation which is adopted for this study because this is exactly the class of Lie groupoids that are relevant in the study of noncommutative geometry of unital $G$-algebras.  

\3 \noindent \textbf{Acknowledgment.} I wish to thank A. Kupiainen and S. Majid for providing me an opportunity to produce this final part of my project.

\section{Noncommutative Orbifolds}

\noindent \textbf{1.1.} It will be assumed that the reader is familiar with the basic concepts of spectral triple and Fredholm module theory, \cite{Con94}. The complex Hilbert spaces associated with these structures are separable and infinite dimensional. Let $G$ be a finite group which acts on the unital algebra $A$. By a noncommutative orbifold we mean a finitely summable spectral triple over the crossed product algebra $G \ltimes A$ with the additional assumption that the $G$-action on the Hilbert space commutes with the Dirac operator of the spectral triple. If a spectral triple has a chiral grading operator $\omega$, then we shall also assume that the $G$-action on the Hilbert space commutes with $\omega$. If a spectral triple does not have a chiral grading, then we write $\omega = 1$. The representation of $G \ltimes A$ on the Hilbert space determines the $G$-module structure of the Hilbert space since we can map $G$ to $G \ltimes A$ by $g \mapsto (g, 1)$ for all $g \in G$. Since the $G$-action commutes with the Dirac operator and $\omega$, these operators can be restricted to the $G$-invariant component of the Hilbert space $\hil^G$, which results in a spectral triple over the invariant subalgebra $A^G$. In particular we have the following structures 
\begin{eqnarray*}
(G \ltimes A, \hil, D, \omega) \hookleftarrow (A^G, \hil^G, D, \omega)
\end{eqnarray*}
where the arrow indicates that $A^G$ and $\hil^G$ can be embedded to $G \ltimes A$ and to $\hil$ as a subalgebra and as a subspace whereas $D$ and $\omega$ commute with the inclusion map of $\hil^G$. The following two cases provide a rich source of noncommutative orbifolds.

\3 \noindent \textbf{Example 1.} Suppose that $G \ltimes X \rightrightarrows X$ is a compact action orbifold (see the notation paragraph in the introduction) which is $G$-spin. Then there is a complex $G$-equivariant Dirac bundle $F_{\Sigma}$ on $X$, and a $G$-invariant Dirac operator $\eth$ acting on the smooth sections of $F_{\Sigma}$, \cite{LM89}. The Dirac operator has a completion to an unbounded and densely defined self-adjoint operator on the Hilbert space of spinors $L^2(F_{\Sigma})$, and there is the following crossed product spectral triple associated with this data: 
\begin{eqnarray*}
(G \ltimes C^{\infty}(X), L^2(F_{\Sigma}), \eth, \omega)
\end{eqnarray*}
If the dimension of $X$ is even, then $L^2(F_{\Sigma})$ is chirally graded. This spectral triple was studied with details in \cite{Har14a} in the more general context of proper {\' e}tale Lie groupoids.  

\3 \noindent \textbf{Example 2.} Let $G$ be a compact simple and simply connected Lie group. The algebras of regular functions on the associated quantum groups $\mathbb{C}[G_q]$ accept a variety of actions under the subgroups of the maximal torus of $G$, \cite{BF12}, \cite{Har14d}. In \cite{Har14d} these algebras were called quantum orbifolds. Suppose that $K$ is a finite subgroup in a maximal torus of $G$ which is given an action on $\mathbb{C}[G_q]$. If this quantum orbifold is spin, $K$ acts on the quantum group effectively, i.e. the $G$ action on the algebra of regular functions $\mathbb{C}[G_q]$ is faithful, then there is a spectral triple on the crossed product algebra $G \ltimes \mathbb{C}[G_q]$ which makes the quantum orbifold a noncommutative orbifold.

\3 In both cases we have assumed an effectiveness condition. In the spectral triple language the effectiveness of a group action on a manifold translates to the faithfulness of the representation of the crossed product algebra on the Hilbert space of spinors, \cite{Har14b}. However, it should be noted that in both cases, Example 1 and Example 2, all the spectral triple axioms apart from the faithfulness of the representation remain to hold if the effectiveness is no longer required.

\3\noindent \textbf{1.2.} The pre-$C^*$-algebras $A$ and $B$ are defined to be Morita equivalent if there is an imprimitivity $B$-$A$ pre-$C^*$-algebra bimodule between them. The definition of such a bimodule is recalled in 2.1. An imprimitivity pre-$C^*$-algebra bimodule has always a completion to an imprimitivity  $C^*$-algebra bimodule between the completed $C^*$-algebras, i.e. the completed $C^*$-algebras are Morita equivalent, \cite{RW98}. If $\hil$ is a Hilbert space representation of the pre-$C^*$-algebra $A$ and if $B$ is Morita equivalent to $A$, i.e. if there is an imprimitivity $B$-$A$ pre-$C^*$-algebra bimodule $E$, then there is an induced map of representations $\hil \mapsto \phi_{\#} \hil$ introduced in \cite{Rie74}. The Hilbert space $\phi_{\#} \hil = E \otimes_A \hil$ is a representation space for $B$ which is defined as follows. One takes the algebraic tensor product $E \odot_A \hil$ over $A$ which accepts a left $B$-module structure since $E$ is a left $B$-module. There is an inner product in $E \odot_A \hil$ which is given by  
\begin{eqnarray*}
\la u \odot \psi, u' \odot \psi' \ra = \la (u',u)_A \psi, \psi' \ra_{\hil}
\end{eqnarray*} 
for all $u,u' \in E$ and $\psi, \psi' \in \hil$, and extended by linearity. Here $(\cdot, \cdot)_A$ denotes the $A$-valued pairing in the bimodule $E$. The completion with respect to this inner product is $\phi_{\#} \hil = E \otimes_A \hil$. 

If $(A, \hil, F)$ is a Fredholm module so that $A,B$ are Morita equivalent pre-$C^*$-algebras through the imprimitivity pre-$C^*$-algebra bimodule $E$, and if $(B, \phi_{\#}\hil, \widetilde{F})$ is a Fredholm module, then $\widetilde{F}$ is called an $F$-connection if the operators
\begin{eqnarray*}
T_u \circ F - \widetilde{F} \circ T_u \5 \text{and} \5 T_u^* \circ \widetilde{F} - F \circ T_u^*
\end{eqnarray*}
are compact for all $u \in E$, where $T_u: \hil \rightarrow E \otimes_A \hil$ is the linear map $\psi \mapsto u \otimes \psi$, and the adjoint of $T_u$ is defined by $T_u^*(v\otimes \psi) = ( u,v )_{A} \psi$, \cite{Ska84}. 

Consider the case of a pair of noncommutative orbifolds associated with the spectral triples 
\begin{eqnarray*}
(G_i \ltimes A_i, \hil_i, D_i, \omega_i) \5 i = 1,2
\end{eqnarray*}
on the crossed product pre-$C^*$-algebras. The spectral triples are defined to be Morita equivalent if the following conditions hold. 

\begin{quote}
\textbf{M1.} There is an imprimitivity $G_2 \ltimes A_2$ - $G_1 \ltimes A_1$ pre-$C^*$-algebra bimodule $E$.

\textbf{M2.} There is a unitary equivalence of spectral triples 
\begin{eqnarray*}
(G_2 \ltimes A_2, \hil_2, D_2, \omega_2) \mapsto (G_2 \ltimes A_2, E \otimes_{G_1 \ltimes A_1} \hil_1, \widetilde{D}_2, \widetilde{\omega}_2),
\end{eqnarray*}
so that $E \otimes_{G_1 \ltimes A_1} \hil_1$ is a $G_2 \ltimes A_2$ representation under the standard left algebra action. 
 
\textbf{M3.} The approximate sign $\widetilde{F}_2$ is an $F_1$-connection on $E \otimes_{G_1 \ltimes A_1} \hil_1$.

\textbf{M4.} The crossed product spectral triples have the same dimension. 

\textbf{M5.} The invariant spectral triples are unitarily equivalent.
\end{quote}
The approximate sign operators in the definition are the bounded Fredholm operators
\begin{eqnarray*}
\widetilde{F}_2 = \frac{\widetilde{D}_2}{(1 + \widetilde{D}_2)^{\frac{1}{2}}} \5 \text{and} \5 F_1 = \frac{D_1}{( 1 + D_1^2)^{\frac{1}{2}}}.
\end{eqnarray*}

\noindent \textbf{Proposition 1.} The Morita equivalence of spectral triples on noncommutative orbifolds is an equivalence relation. 

\3 \noindent Proof. The property of having an imprimitivity pre-$C^*$-algebra bimodule defines an equivalence relation, \cite{RW98} Proposition 3.16. In particular, for the reflexivity we can use the algebra $G \ltimes A$ itself as a pre-$C^*$-imprimitivity bimodule. Then \textbf{M2} holds with $\hil \simeq G \ltimes A \otimes_{G \ltimes A} \hil$ given by  $\psi \mapsto \textbf{1} \otimes \psi$, and \textbf{M3} holds because the commutators $[F, a]$ are compact operators for all $a \in G \ltimes A$. Therefore the axioms \textbf{M1-3} are reflexive. Suppose then that the  imprimitivity bimodule $E$ has been fixed. Then there are the crossed product $C^*$-algebras, $G_1 \ltimes \mathcal{A}_1$ and $G_2 \ltimes \mathcal{A}_2$, and the completion  $\mathcal{E}$ of $E$ is a  $G_2 \ltimes \mathcal{A}_2$ - $G_1 \ltimes \mathcal{A}_1$ $C^*$-algebra bimodule. The bimodule defines a class $(\mathcal{E},G_2 \ltimes \mathcal{A}_2, G_1 \ltimes \mathcal{A}_1,0)$ in the operator KK-theory group $KK_0(G_2 \ltimes \mathcal{A}_2,G_1 \ltimes \mathcal{A}_1)$, where the operator is $0$. By definition, the conditions $\textbf{M2}$ and $\textbf{M3}$ mean that, up to a unitary equivalence, $(G_2 \ltimes \mathcal{A}_2, \hil_2, F_2)$ is a representative of the intersection product of $(\mathcal{E},G_2 \ltimes \mathcal{A}_2, G_1 \ltimes \mathcal{A}_1,0)$ with $(G_1 \ltimes \mathcal{A}_1, \hil_1, F_1)$. This property is transitive with respect to successive Morita bimodule operations. In addition, the intersection product under the class of the dual Morita bimodule of $\mathcal{E}$ inverts the relation above. Therefore, \textbf{M1-3} define an equivalence relation. 

Clearly we can strengthen this equivalence relation with \textbf{M4-5} to get a new equivalence relation. \5 $\square$

\3 The condition \textbf{M1} implies that the usual homotopy invariants of algebras, such as the cyclic homology of the pre-$C^*$-algebra, or the $K$-theory of the $C^*$-algebraic completion are invariant under the Morita equivalence of spectral triples on noncommutative orbifolds. The conditions \textbf{M2-3} imply that the K-homology class associated with the crossed product spectral triple is a Morita invariant. Condition \textbf{M5} implies that the invariant spectral geometry is invariant. All these properties are very realistic from the viewpoint of topology and differential geometry of orbifolds: topological homology invariants, topological K-theory and K-homology are known to be Morita invariant, and so is the orbit space as a metric space which is modelled with the invariant spectral triple. One might argue that the dimensionality condition \textbf{M4} is too strong. Namely, the dimensionality of a groupoid base manifolds is not a Morita invariant in general. However, in the case under consideration, the geometric theory is restricted to global action groupoids subject to a finite group action, and in this case all Morita equivalences are {\' e}tale structure preserving, and the dimensionality is indeed invariant. So \textbf{M4} should be implemented in the case under consideration.

\3 \noindent \textbf{1.3.} Suppose that $G \ltimes X$ and $K \ltimes Y$ is a pair of Morita equivalent compact action orbifolds and $G \ltimes X$ is $G$-spin. There is a Dirac spectral triple on $G \ltimes C^{\infty}(X)$, as in Example 1. Since the $G$-spin property can be described in terms of groupoid hypercohomology, it follows that $K \ltimes Y$ is $K$-spin, \cite{Har14a}. There is a geometric Morita equivalence (a Morita bitorsor) between these groupoid. The definition will be recalled in 2.2.  This geometric Morita equivalence, to be denoted by $\phi$, allows one to define a $K$-equivariant induced spinor bundle $\phi_{\#} F_{\Sigma}$ on $Y$, see 2.4. Similarly, there are the induced $K$-invariant Dirac operator $\phi_{\#} \eth$ and the induced $K$-invariant chirality grading operator $\phi_{\#} \omega$ on $\phi_{\#} F_{\Sigma}$, see 3.1. In this notation, we state the main result of this manuscript: 

\3 \noindent \textbf{Theorem 1.} The Dirac spectral triples 
\begin{eqnarray*}
(G \ltimes C^{\infty}(X), L^2(F_{\Sigma}), \eth, \omega) \5 \text{and}\5 (K \ltimes C^{\infty}(Y), L^2(\phi_{\#} F_{\Sigma}), \phi_{\#}\eth, \phi_{\#}\omega)
\end{eqnarray*}
are Morita equivalent as noncommutative orbifolds.  

\3 A proof for this theorem occupies the sections 2 and 3 of this manuscript, and is put together in the subsection 3.2.  

\3 \noindent \textbf{1.4.} Let us then discuss the case where $X/G$ is a smooth manifold, or equivalently, the case where $G$ acts freely on $X$. These conditions are equivalent to the existence of a geometric Morita equivalence between the action groupoids $G \ltimes X$ and $1 \ltimes X/G$. The latter groupoid is the unit groupoid of the manifold $X/G$, i.e. there are only unit arrows. In particular we can view $G \ltimes X$ as a representative of a manifold object in the category of Lie groupoids. 

In the spectral triple language of 1.2, we consider a noncommutative orbifold represented by a spectral triple $(G \ltimes A, \hil, \eth, \omega)$ to be a smooth noncommutative geometric space, if this spectral triple is Morita equivalent to the invariant spectral triple $(A^G, \hil^G, \eth, \omega)$ (identify $1 \ltimes A^G$ with $A^G$). This is how we can reduce the equivariant crossed product formalism to the usual spectral triple formalism where the $G$-action is absent. Notice that in this case the equivariant homotopy invariants, such as the cyclic homology or the K-theory invariants associated with the crossed product algebra $G \ltimes A$, are isomorphic to the corresponding invariants computed from the algebra $A^G$. Therefore, exactly as in the topological case, the equivariant formalism becomes redundant as soon as this freeness condition is satisfied. 

Smoothness of noncommutative orbifolds has been also studied in the recent manuscript \cite{BS14} from an algebraic geometric point of view. It would be interesting to understand to what extent this theory coincides with the differential geometric approach introduced above.

\section{Induced Representations}

\noindent \textbf{2.1.} Let $A$ and $B$ be pre-$C^*$-algebras. The $C^*$-algebraic completions are denoted by $\mathcal{A}$ and $\mathcal{B}$. A $B$-$A$-bimodule $E$ is an imprimitivity pre-$C^*$-algebra bimodule if the following axioms hold  
\begin{quote}
\textbf{1.} $E$ is a right pre-inner product $A$-module and a left pre-inner product $B$-module: there are the pairings $(\cdot, \cdot)_A : E \times E \rightarrow A$  and ${}_B(\cdot, \cdot) : E \times E \rightarrow B$ so that
\begin{eqnarray*}
(u, \lambda v + \mu w)_A = \lambda (u,v)_A + \mu (u,w)_A && {}_B(\lambda u + \mu v, w) = \lambda {}_B(u,w) + \mu {}_B(v,w), \\
(u, va )_A = (u,v)_A a  && {}_B(b u, v) = b{}_B(u,v) \\
(u,v)^*_A = (v,u)_A && {}_B(u,v)^* = {}_B(v,u)  \\
(u,u)_A \geq 0 && {}_B(u,u) \geq 0  
\end{eqnarray*}
hold for all $\lambda,\mu \in \mathbb{C}$, $a \in A$, $b \in B$ and $u,v,w \in E$, and the inequalities hold in the $C^*$-algebra completions $\mathcal{A}$ and $\mathcal{B}$. 

\textbf{2.} The linear spans ${}_B(u,v)$ and $(u,v)_A$ with $u,v \in E$ are dense in $\mathcal{A}$ and in $\mathcal{B}$. 

\textbf{3.} For all $a \in A$, $b \in B$ and $u \in E$: 
\begin{eqnarray*}
(b u, b u)_A \leq ||b||^2(u,u)_A \5 \text{and} \5 {}_B(ua,ua) \leq ||a||^2{}_B(u,u).
\end{eqnarray*}

\textbf{4.} For all $u,v,w \in E$: $u (v,w)_A = {}_B (u,v) w $. 
\end{quote}

\noindent One can apply the $C^*$-algebra norms in $\mathcal{A}$ and in $\mathcal{B}$ with the pairings ${}_A(\cdot, \cdot)$ and $(\cdot, \cdot)_B$ to obtain a pair of norms in $E$. These norms are equivalent. If $E$ is a pre-imprimitivity bimodule, then it has a completion to an imprimitivity $C^*$-algebra bimodule $\mathcal{E}$, see \cite{RW98} Proposition 3.12. In particular, the $C^*$-algebras $\mathcal{A}$ and $\mathcal{B}$ are Morita equivalent.  

\3 \noindent \textbf{2.2.} Let $G \ltimes X$ and $K \ltimes Y$ be a pair of compact action orbifolds. These orbifolds are defined to be Morita equivalent if there is a Morita bitorsor, i.e. a smooth manifold $Q$ and a diagram
\begin{center}
\gpdhom{Q}{K \ltimes Y}{Y}{X}{G \ltimes X}{\alpha}{\varrho}
\end{center}
so that the following conditions are satisfies
\begin{quote}
\textbf{1.} The action groupoid $G \ltimes X$ acts on $Q$ from the right such that the anchor for this action is $\varrho$. The action groupoid $K \ltimes Y$ acts on $Q$ from the left such that the anchor for this action is $\alpha$. The actions are mutually commutative. 

\textbf{2.} The maps $\varrho$ and $\alpha$ are local diffeomoprhisms, the action of $G \ltimes X$ is free and transitive on the fibres $\alpha^{-1}(y)$ for all $y\in Y$, and the action of $K \ltimes Y$ is free and transitive on the fibres $\varrho^{-1}(x)$ for all $x \in X$. 
\end{quote}

\noindent Let us clarify these conditions slightly. The domain for the left action on $Q$ is the fibre product
\begin{eqnarray*}
(K \times Y) \fp{s}{\alpha} Q
\end{eqnarray*}
So, if $q \in Q$ is fixed, then the action restricts to the map 
\begin{eqnarray*}
K \times \{q\} \rightarrow \varrho^{-1}(\varrho(q)); \5 (k, q) \mapsto (k, \alpha(q)) \cdot q
\end{eqnarray*}
 This map is well defined since $s(k, \alpha(q)) = \alpha(q)$ for all $k \in K$, and it follows from \textbf{2} that the image is the set $\varrho^{-1}(\varrho(q))$. Since the action is free and transitive, these maps are bijective for all $\{q\}$, and all the $\varrho$-fibres in $Q$ have the cardinality equal to $\# K$, i.e. the number of group elements in $K$. Observe that if $\tau: y \mapsto \tau \cdot y$ is an arrow and if $\alpha(q) = y$, then $\alpha(\tau \cdot q) = \tau \cdot \alpha(q)$. Similarly, the domain of the right action is given by 
\begin{eqnarray*}
(G \times X) \fp{t}{\varrho} Q
\end{eqnarray*}
If $q \in Q$, then the right action restricts to define the map
\begin{eqnarray*}
\{q\} \times G \rightarrow \alpha^{-1}(\alpha(q)); \5 (q,g) \mapsto q \cdot (\sigma, \sigma^{-1} \cdot \varrho(q))
\end{eqnarray*}
This map is well defiend since $t(\sigma, \sigma^{-1} \cdot \varrho(q)) = \varrho(q)$ for all $\sigma \in G$, and by \textbf{2} this is a bijective map onto $\alpha^{-1}(\alpha(q))$. The cardinality of an arbitrary $\alpha$-fibre is equal to $\# G$. Observe that if $\sigma^{-1}: x \rightarrow \sigma^{-1} \cdot x$ is an arrow, and if $\varrho(q) = x$, then $\varrho(q \cdot \sigma) = \sigma^{-1} \cdot \varrho(q)$. 

The Morita equivalence has a particularly convenient geometric structure in the case of compact action orbifolds. 

\3 \noindent \textbf{Proposition 2.} Suppose that $G \ltimes X$ and $K \ltimes Y$ are finite action orbifolds. If $Q$ is a Morita bitorsor between $G \ltimes X$ and $K \ltimes Y$, then $\varrho: Q \rightarrow X$ is an $\# K$-sheeted covering space, and $\alpha: Q \rightarrow Y$ is a $\# G$-sheeted covering space and $Q$ is compact.

\3 \noindent Proof. According to the analysis above, the cardinality of $\varrho^{-1}(x)$ is constant for all $x \in X$, and equal to $\# K$. It follows that $\varrho$ is a covering projection, \cite{Ho75} Lemma 2. Similarly, we have seen that $\alpha^{-1}(y)$ is constant for all $y \in Y$, and equal to $\# G$, and so $\alpha$ is a covering projection. $Q$ is necessarily compact since it is a smooth covering space with a finite number of sheets over a compact manifold. \5 $\square$ 

\3 \noindent \textbf{2.3.} In what follows we need to apply Haar integrals over the subsets such as $(G \ltimes X)^x$ and $(K \ltimes Y)^y$ for some $x \in X$ and $y \in Y$, and for the sake of simplicity of the notation, we shall write
\begin{eqnarray*}
\Theta = G \ltimes X \5 \text{and} \5 \Xi = K \ltimes Y. 
\end{eqnarray*}
Suppose that $Q$ is a Morita bitorsor, as in 2.2. Then we equip $C^{\infty}(Q)$ with the left $C^{\infty}(\Xi)$ and the right $C^{\infty}(\Theta)$ module structures which are determined by  
\begin{eqnarray}\label{bimodule}
(f \cdot a)(q) &=& \int_{\Theta^{\varrho(q)}} a(\sigma^{-1}) f(q \cdot \sigma), \nonumber \\
(b \cdot f)(q) &=& \int_{\Xi^{\alpha(q)}} b(\tau) f(\tau^{-1} \cdot q)  
\end{eqnarray}
for all $f \in C^{\infty}(Q)$, $a \in C^{\infty}(\Theta)$, $b \in C^{\infty}(\Xi)$ and $q \in Q$. There is the $C^{\infty}(\Theta)$-valued pairing in $C^{\infty}(Q)$ defined by
\begin{eqnarray*}
( f, g )_{\Theta}(\sigma) = \int_{\Xi^{\alpha(q)}} \overline{f(\tau^{-1} \cdot q)} g(\tau^{-1} \cdot q \cdot \sigma)
\end{eqnarray*}
where $q$ is any point so that $t(\sigma) = \varrho(q)$. Since the $\Xi$-action is transitive on the set of points with this property, the choice is arbitrary. There is also the $C^{\infty}(\Xi)$-valued pairing in $C^{\infty}(Q)$ defined by
\begin{eqnarray*}
{}_{\Xi}(f,g)(\tau) = \int_{\Theta^{\varrho(q)}} f(\tau^{-1} \cdot q \cdot \sigma) \overline{g(q \cdot \sigma)}
\end{eqnarray*}
where $q$ is any point so that $t(\tau) = \alpha(q)$, and the choice is arbitrary. These structures were introduced in \cite{MRW89}. 

\3 \noindent \textbf{Proposition 3.} If $G \ltimes X$ and $K \ltimes Y$ are Morita equivalent compact action orbifolds, then $C^{\infty}(Q)$ equipped with the bimodule structure \eqref{bimodule} and pairings $(\cdot, \cdot)_{\Theta}$ and ${}_{\Xi}(\cdot,\cdot)$ is an imprimitivity pre-$C^*$-algebra $C^{\infty}(K \ltimes Y)$ - $C^{\infty}(G \ltimes X)$-bimodule. 

\3 \noindent Proof. The corresponding claim for proper topological groupoids has been proved in the reference \cite{MRW89}. So, in the smooth case under consideration, the axioms \textbf{1}, \textbf{3} and \textbf{4} of 2.1 hold because they are valid for all continuous functions.  

In the reference \cite{MRW89}, the algebra $C(Q)$ of continuous functions in $Q$ is equipped with the inductive limit topology. However, by Proposition 2 the covering space $Q$ is compact and therefore the inductive limit topology in $C(Q)$ is the usual norm topology. In particular, the algebra of smooth functions $C^{\infty}(Q)$ equipped with the standard Frechet topology is dense in $C(Q)$. Therefore the axiom \textbf{2} of 2.1 is fulfilled in the smooth case as well. $\square$

\3 \noindent \textbf{2.4.} Suppose that $\xi$ is a complex vector bundle on $G \ltimes X$. More precisely $\pi: \xi \rightarrow X$ is a complex vector bundle and $G \ltimes X$ acts on $\xi$ with respect to the anchor map $\pi$. These bundles can be identified with $G$-equivariant vector bundles. If $\sigma \in \Theta_x$, then let us denote by $\rho(\sigma)$ the linear isomorphism $\xi_x \rightarrow \xi_{\sigma \cdot x}$ associated with the action. Denote by $\Gamma^{\infty}(\xi)$ the space of smooth sections in $\xi$. The space of smooth sections is a $C^{\infty}(G \ltimes X)$ module under the action
\begin{eqnarray*}
(a \cdot \psi)_x = \int_{\Theta^x} a(\sigma) (\varphi^{\#}_{\sigma^{-1}} \psi)_x
\end{eqnarray*}
where we are writing $(\varphi^{\#}_{\sigma^{-1}} \psi)_x = \rho(\sigma) \psi_{\sigma^{-1} \cdot x}$ for all $x \in \Theta^x$. Observe that this is an element in the fibre of $\xi$ at $x \in X$ because $\rho(\sigma): \xi_{\sigma^{-1} \cdot x} \rightarrow \xi_x$ is a linear isomorphism. 

Suppose that $G \ltimes X$ and $K \ltimes Y$ are Morita equivalent compact action orbifolds. Let us fix a Morita equivalence $\phi = (\alpha, Q, \varrho)$ as in 2.2. This Morita equivalence can be used to push bundles on $G \ltimes X$ to bundles on $K \ltimes Y$. If $\xi$ is a  complex vector bundle on $G \ltimes X$, then define 
\begin{eqnarray*}
\phi_{\#} \xi = [Q \fp{\varrho}{\pi} \xi]/G \ltimes X. 
\end{eqnarray*}
which is a vector bundle on $K \ltimes Y$ of equal rank. The $G \ltimes X$-action is defined by 
\begin{eqnarray*}
\sigma \cdot [q, u] = [q \cdot \sigma^{-1}, \rho(q)u]
\end{eqnarray*}
 for all $\sigma \in (G \ltimes X)_{\varrho(q)}$ and $u \in \xi_{\varrho(q)}$. The bundle projection $\pi': \phi_{\#} \xi \rightarrow Y$ is the map $\pi': [q, u] \mapsto \alpha(q)$ and the left $\Xi$ action is given by $\tau \cdot [q, u] = [\tau \cdot q, u]$ for all $\tau \in (K \ltimes Y)_{\alpha(q)}$. The fibre of $\phi_{\#} \xi$ at $y\in Y$ is an equivalence class of fibres of $\xi$ at the points $\varrho(q)$ for all $q \in Q$ so that $q \in \alpha^{-1}(y)$. There are exactly $\# G$ such points by Proposition 2. If $\eta$ is a section of $\phi_{\#} \xi$, then we write $\eta_q(y)$ for the value of $\eta$ at $y$ in the representative $\xi_{\varrho(q)}$ for the fibre. In this notation, the sections satisfy
\begin{eqnarray*}
\rho(\sigma) \eta_{q \cdot \sigma}(y) = \eta_q(y)
\end{eqnarray*}
for all $\sigma \in \Theta^{\varrho(q)}$. 

Let us denote by $\phi_{\#} \Gamma^{\infty}(\xi) = C^{\infty}(Q) \odot_{\Theta} \Gamma^{\infty}(\xi)$ the algebraic tensor product of $C^{\infty}(Q)$ and $\Gamma^{\infty}(\xi)$ over the algebra $C^{\infty}(G \ltimes X)$. 
 
\3 \noindent \textbf{Proposition 4.} There is a $\mathbb{C}$-linear map $\chi: \phi_{\#} \Gamma^{\infty}(\xi) \rightarrow  \Gamma^{\infty}(\phi_{\#} \xi)$ which sends the class of $f \odot \psi$ in $\phi_{\#} \Gamma^{\infty}(\xi)$ to the section of $\phi_{\#} \xi$ which is defined by
\begin{eqnarray*}
\chi(f \odot \psi)_q(y) = \int_{\Theta^{\varrho(q)}} f(q \cdot \sigma) (\varphi^{\#}_{\sigma^{-1}} \psi)_{\varrho(q)} = \frac{1}{\# G} \sum_{\sigma \in \Theta^{\varrho(q)}} f(q \cdot \sigma) (\varphi^{\#}_{\sigma^{-1}} \psi)_{\varrho(q)} 
\end{eqnarray*}
for all $y \in Y$ and $q \in \alpha^{-1}(y)$. \3

\noindent Proof. The integration and the pullbacks are linear operations and therefore $\chi$ is a linear map. If $f \odot \psi$ is a representative of a class in $\phi_{\#} \Gamma^{\infty}(\xi)$, then 
\begin{eqnarray*}
\rho(\sigma) \chi(f \odot \psi)_{q \cdot \sigma}(y) &=& \rho(\sigma) \int_{\tau \in \Theta^{\varrho(q \cdot \sigma)}} f(q \cdot \sigma \cdot \tau) (\varphi^{\#}_{\tau^{-1}} \psi)_{\varrho(q \cdot \sigma)} \\
 &=& \rho(\sigma) \int_{\tau \in \Theta^{\varrho(q \cdot \sigma)}} f(q \cdot \sigma \cdot \tau) (\varphi^{\#}_{\tau^{-1}} \psi)_{\sigma^{-1} \cdot \varrho(q)} \\
  &=& \int_{\tau \in \Theta^{\varrho(q \cdot \sigma)}} f(q \cdot \sigma \cdot \tau) \rho(\sigma) \rho(\tau) (\varphi^{*}_{(\sigma \cdot \tau)^{-1}} \psi)_{\varrho(q)} \\
	&=&  \int_{\tau \in \Theta^{\varrho(q \cdot \sigma)}} f(q \cdot \sigma \cdot \tau) (\varphi^{\#}_{(\sigma \cdot \tau)^{-1}} \psi)_{\varrho(q)} \\
	&=& \int_{\tau \in \Theta^{\varrho(q)}} f(q \cdot \tau) (\varphi^{\#}_{\tau^{-1}} \psi)_{\varrho(q)} \\
  &=& \chi(f \odot \psi)_q(y).
\end{eqnarray*}
for all $q \in  \alpha^{-1}(y)$ and $\sigma \in \Theta^{\varrho(q)}$. It follows by linearity that the image of $\chi$ lies in $\Gamma^{\infty}(\phi_{\#} \xi)$. 

It remains to show that the image of $\chi$ is independent on the choices of the equivalence classes in $\phi_{\#} \Gamma^{\infty}(\xi)$. Let $f \odot \psi$ be a represenative of an element in $C^{\infty}(Q) \odot_{\Theta} \Gamma^{\infty}(\xi)$ and $a \in C^{\infty}(G \ltimes X)$. Then $f \odot a \cdot \psi$ maps to the following element under $\chi$:
\begin{eqnarray*}
\chi(f \odot a \cdot \psi)_q(y) &=& \int_{\sigma \in \Theta^{\varrho(q)}} f(q \cdot \sigma) \int_{\tau \in \Theta^{\sigma^{-1} \cdot x}} a(\tau) (\varphi^{\#}_{\sigma^{-1}}\varphi^{\#}_{\tau^{-1}} \psi)_{\varrho(q)} \\
&=& \int_{\sigma \in \Theta^{\varrho(q)}} \int_{\tau \in \Theta^{\varrho(q \cdot \sigma)}} f(q \cdot \sigma)  a(\tau) (\varphi^{\#}_{(\sigma \cdot \tau)^{-1}}\psi)_{\varrho(q)}
\end{eqnarray*}
and $f \cdot a \odot \psi$ maps to 
\begin{eqnarray*}
\chi(f \cdot a \odot \psi)_q(y) &=& \int_{\sigma \in \Theta^{\varrho(q)}} \int_{\tau \in \Theta^{\sigma^{-1} \cdot \varrho(q)}} f(q \cdot \sigma \cdot \tau) a(\tau^{-1}) (\varphi^{\#}_{\sigma^{-1}} \psi)_{\varrho(q)}\\
&=& \int_{\sigma \in \Theta^{\varrho(q)}} \int_{\tau \in \Theta^{\varrho(q \cdot \sigma)}} f(q \cdot \sigma \cdot \tau) a(\tau^{-1}) (\varphi^{\#}_{\sigma^{-1}} \psi)_{\varrho(q)}
\end{eqnarray*}
In the first case, the argument of $f$ runs over $q \cdot \sigma$, the argument of $a$ runs over the arrows $\varrho(q \cdot \sigma \cdot \tau) \rightarrow \varrho(q \cdot \sigma) $, and the argument of $\varphi^{\#}$ runs over the arrows $\varrho(q \cdot \sigma \cdot \tau) \rightarrow \varrho(q)$.  In the second case, the argument of $f$ runs over $q \cdot \sigma \cdot \tau$, the argument of $a$ runs over $\varrho(q \cdot \sigma) \rightarrow \varrho(q \cdot \sigma \cdot \tau) $ and $\varphi^{\#}$ runs over $\varrho(q \cdot \sigma) \rightarrow \varrho(q)$. The claim follows since the fibre integration is with respect to the counting measure in $G$. \5 $\square$

\3 Suppose that $\xi$ is a complex vector bundle on $G \ltimes X$ which is equipped with a smoothly varying $G$-invariant hermitian inner product $(\cdot, \cdot)$ in its fibres. This gives the pairing $\Gamma^{\infty}(\xi) \otimes \Gamma^{\infty}(\xi) \rightarrow C^{\infty}(X)$. If $\phi = (\alpha, Q, \varrho)$ is a Morita equivalence as above, then there is a hermitian structure $(\cdot, \cdot)_{\#}$ in the fibres of $\phi_{\#} \xi$ which is defined by
\begin{eqnarray*}
([v_1], [v_2])_{\#, y} = (v^{\varrho(q)}_1, v_2^{\varrho(q)})_{\varrho(q)}
\end{eqnarray*}
where $q$ is any point in $\alpha^{-1}(y)$ and $v_i^{\varrho(q)}$ is the class of $[v_i]$, $i = 1,2$ evaluated at the fibre $\xi_{\varrho(q)}$. The value of the inner product is independent on the choice of such $q$ because of the $G$-invariance. This gives us the pairing $\Gamma^{\infty}(\phi_{\#} \xi) \otimes \Gamma^{\infty}(\phi_{\#} \xi) \rightarrow C^{\infty}(Y)$. We also assume that $G \ltimes X$ is a riemannian groupoid, and that $\nu$ is a $G$-invariant riemannian volume form. Then there is the induced $K$-invariant riemannian volume form $\phi_{\#} \nu$ on $Y$. This is the differential form on $Y$ which is uniquely determined by the equality $\varrho^*(\nu) = \alpha^*(\phi_{\#} \nu)$. One can apply local sections of $\alpha$ to pull $\varrho^*(\nu)$ to a globally defined form on $Y$, \cite{Har14b}. 

\3 We define the Hilbert space $\phi_{\#} L^2(\xi)$ to be the completion of $C^{\infty}(Q) \odot_{\Theta} \Gamma^{\infty}(\xi)$ in the norm determined by the inner product 
\begin{eqnarray*}
\la f_1 \odot \psi_1, f_2 \odot \psi_2 \ra &=& \la(f_2, f_1)_{\Theta} \cdot \psi_1, \psi_2 \ra_X \\
&=& \int_X ((f_2, f_1)_{\Theta} \cdot \psi_1, \psi_2) \nu.  
\end{eqnarray*} 
The bounded left action of $C^{\infty}(\Xi)$ on  $C^{\infty}(Q) \odot_{\Theta} \Gamma^{\infty}(\xi)$ extends to make $\phi_{\#} L^2(\xi)$ a representation space for $C^{\infty}(\Xi)$. There is also a natural inner product in $\Gamma^{\infty}(\phi_{\#} \xi)$ given by
\begin{eqnarray*}
\la \eta_1, \eta_2 \ra_Y = \int_Y (\eta_1, \eta_2)_{\#} (\phi_{\#}\nu). 
\end{eqnarray*}
for all $\eta_1, \eta_2 \in \Gamma^{\infty}(\phi_{\#} \xi)$. The Hilbert space completion of $\Gamma^{\infty}(\phi_{\#} \xi)$ in the norm determined by this inner product is denoted by $L^2(\phi_{\#} \xi)$. 

\3 \noindent \textbf{Theorem 2.} Let $G \ltimes X$ and $K \ltimes Y$ be Morita equivalent compact action orbifolds. The map $\chi$ extends to an isomorphism of $C^{\infty}(K \ltimes Y)$-representations: 
\begin{eqnarray*}
\chi: \phi_{\#} L^2(\xi) \rightarrow  L^2(\phi_{\#} \xi)
\end{eqnarray*}
for any vector bundle $\xi$ on $G \ltimes X$ and Morita equivalence $\phi$. 

\3 \noindent Proof. We shall prove the injectivity of $\chi: \phi_{\#} \Gamma^{\infty}(\xi) \rightarrow  \Gamma^{\infty}(\phi_{\#} \xi)$ first. Let $f \odot \psi$ be a representative of a class in $C^{\infty}(Q) \odot_{\Theta} \Gamma^{\infty}(\xi)$. Using the $C^{\infty}(G \ltimes X)$-linearity of the tensor product, we find the following equalities
\begin{eqnarray*}
f(q) \odot \psi_x &=& \frac{1}{\# G} \sum_{\sigma \in \Theta^{\varrho(q)}} f(q) \cdot \textbf{1}_{\sigma^{-1}} \cdot \textbf{1}_{\sigma} \odot \psi_x \\
&=& \frac{1}{\# G} \sum_{\sigma \in \Theta^{\varrho(q)}} f(q) \cdot \textbf{1}_{\sigma^{-1}}  \odot \textbf{1}_{\sigma} \cdot \psi_x \\
&=&  \frac{1}{\# G} \sum_{\sigma \in \Theta^{\varrho(q)}} f(q \cdot \sigma)  \odot (\varphi^{\#}_{\sigma^{-1}} \psi)_x 
\end{eqnarray*}
and consequently $\chi(f \odot \psi) = 0$ implies that $f(q) \odot \psi_{\varrho(q)} = 0$ for all $q \in Q$. Using the linearity of $\chi$ we also have that
\begin{eqnarray*}
\chi(\sum_i f^i \odot \psi^i) = 0 \5 \Rightarrow \5 \sum_i f^i(q) \odot \psi^i_{\varrho(q)} = 0.
\end{eqnarray*}
for all $q \in Q$. So, it is sufficient to show that $\sum f^i(q) \odot \psi^i_{\varrho(q)} = 0$ holds only if $\sum f^i \odot \psi^i$ is the zero vector of $C^{\infty}(Q) \odot_{\Theta} \Gamma^{\infty}(\xi)$.

Suppose that $\psi$ is a section of $\xi$ so that its support lies in a compact subset $V$ of $X$ and $V$ is a subset in an open subset $U$ of $X$ where $\xi$ is a trivial bundle. On $V$ we can choose the frame coordinate fields $e_i$, $i \in \{1, \ldots, \text{rank}(\xi)\}$ so that $\{(e_i)_x\}$ span the fibre of $\xi_x$ for all $x \in V$. The frame fields can be extended to smooth vector fields on $X$ so that their supports lie in $U$. Then we can write $\psi = \sum \psi^i e_i$ where $\psi^i$ are smooth functions with supports in $U$, and
\begin{eqnarray*}
f \odot \psi &=& f \odot \sum_{i = 1}^{\text{rank}(\xi)} \psi^i e_i \\
&=& f \odot \sum_{i = 1}^{\text{rank}(\xi)} (\psi^i \circ t) \textbf{1}_{e} \cdot e_i \\
&=& \sum_{i = 1}^{\text{rank}(\xi)} f \cdot (\psi^i \circ t) \textbf{1}_{e} \odot e_i \\
&=& \sum_{i = 1}^{\text{rank}(\xi)} f(\psi^i \circ \varrho) \odot e_i.
\end{eqnarray*}
Thus, if $f(\psi_1^i \circ \varrho) = 0$ as a function on $Q$, or equivalently $f(q) \odot \psi_{\varrho(q)} = 0$ for all $q \in Q$, then $f \odot \psi = 0$. If $\psi$ is not supported in a neighborhood of $X$ where $\xi$ trivializes, then the standard partition of unity argument can be used to generalize the argument to this case straightforwardly. The general case holds by linearity. The injectivity of $\chi: \phi_{\#} \Gamma^{\infty}(\xi) \rightarrow  \Gamma^{\infty}(\phi_{\#} \xi)$ follows.  

The next task is to prove the surjectivity of $\chi: \phi_{\#} \Gamma^{\infty}(\xi) \rightarrow  \Gamma^{\infty}(\phi_{\#} \xi)$. Let $\xi \rightarrow X$ be a vector bundle on $G \ltimes X$. Let us fix an open cover $\{U_a: a \in I\}$ for $X$ and choose the components to be small enough so that $\xi$ can be trivialized over each $U_a$, 
\begin{eqnarray*}
\varrho^{-1}(U_{a}) = \coprod_{i = 1}^{\# K} Q_{a}^i
\end{eqnarray*}
where $Q_a^i$ are mutually disconnected open subsets of $Q$, and $\varrho$ restricts to a diffeomorphism $Q_a^i \rightarrow U_a$ for each $a$ and $i$. The components $Q_a^i$ exist because $Q$ is a $\# K$ sheeted covering space and $\varrho: Q \rightarrow X$ is the covering projection. In addition, by shrinking $U_a$ further if necessary, we can assume that the $G$ orbit of each $Q_a^i$ consist of $\# G$ disconnected components. This assumption holds if $U_a$ are chosen to be small enough because the $G$ action on $Q$ is free and transitive. Now $\{Q_a^i: a \in I, 1 \leq i \leq \# K \}$ is an open cover of $Q$, and $\{\alpha(Q_a^i): a \in I, 1 \leq i \leq \# K \}$ is an open cover of $Y$ because $\alpha$ is a covering map and so an open map. Let us fix a partition of unity $\rho_a^i$ subordinate to the cover $\{Q_a^i\}$ of $Q$.

Suppose that $\eta$ is a section of $\phi_{\#} \xi$. The bundle $\phi_{\#} \xi$ trivializes over each $\alpha(Q_a^i)$ because the fibre $(\phi_{\#} \xi)_y$ at $y \in \alpha(Q_a^i)$ is represented by the fibres of $\xi$ at $\varrho(q) \in U_a$ for $q$ such that $q \in \alpha^{-1}(y)$, and $\xi$ trivializes over $U_a$. So, if $q' \in Q_a^i$, then we have the local neighborhood $\alpha(Q_a^i)$ of $\alpha(q')$ and for all $y \in \alpha(Q_a^i)$ we can write 
\begin{eqnarray*}
\eta_q(y) = (\psi_a^i)_{\varrho(q)}
\end{eqnarray*}
where $\psi_a^i$ is a locally defined section $\psi_i^a: U_a \rightarrow \xi$. We can go through all $Q_a^i$ which results in locally defined sections $\psi_a^i$ for all $a\in I$ and $1 \leq i \leq \# K$. Since $\varrho$ restricts to a diffeomorphism over each $Q_a^i$ there are the unique local sections $\delta_a^i: U_a \rightarrow Q_a^i$ of $\varrho$ for all $a$ and $i$.  Then we define the vector 
\begin{eqnarray*}
\sum_{a,i} \sqrt{\rho_a^i} \odot \sqrt{\rho_a^i \circ \delta_a^i} \psi_a^i \in C^{\infty}(Q) \odot_{\Theta} \Gamma^{\infty}(\xi).  
\end{eqnarray*}
Its image under $\chi$ is given by 
\begin{eqnarray*}
\chi \Big( \sum_{a,i} \sqrt{\rho_a^i} \odot \sqrt{\rho_a^i \circ \delta_a^i} \psi_a^i \Big)_q(y) &=& \int_{\sigma \in \Theta^x} \sum_{a,i} \sqrt{\rho_a^i}(q \cdot \sigma) \varphi_{\sigma^{-1}}^{\#} (\sqrt{\rho_a^i \circ \delta_a^i} \psi_a^i)_{\varrho(q)}\\ 
&=&  \sum_{a,i} \int_{\sigma \in \Theta^x}\sqrt{\rho_a^i}(q \cdot \sigma) \varphi_{\sigma^{-1}}^{\#} (\sqrt{\rho_a^i \circ \delta_a^i} \psi_a^i)_{\varrho(q)} \\
&=& \sum_{a,i}  \sqrt{\rho_a^i}(q ) \sqrt{\rho_a^i \circ \delta_a^i}(\varrho(q)) (\psi_a^i)_{\varrho(q)} \\
&=& \sum_{a,i} \rho_a^i(q) (\psi_a^i)_{\varrho(q)} = \eta_q(y)
\end{eqnarray*}
for all $y \in Y$ and $q \in \alpha^{-1}(y)$. The third equality holds because the cover $Q_a^i$ is chosen so that the $G$-orbits of the components are mutually disjoint, and the support of $\rho_a^i$ lies in $Q_a^i$ for all $a$ and $i$. The surjectivity of $\chi: \phi_{\#} \Gamma^{\infty}(\xi) \rightarrow  \Gamma^{\infty}(\phi_{\#} \xi)$ follows.

The next task is to show that the module structures are equivalent. Let $f \odot \psi$ be a representative of a class in $C^{\infty}(Q) \odot_{\Theta} \Gamma^{\infty}(\xi)$ For any $b \in C^{\infty}(K \ltimes Y)$ we have 
\begin{eqnarray*}
 b \cdot (f(q) \odot \psi_x) =  \int_{\Xi^{\alpha(q)}} b(\tau)f(\tau^{-1} \cdot q) \odot \psi_{x}.
\end{eqnarray*} 
Under the vector space isomorphism $\chi$, this element maps to the section of $\phi_{\#} \xi$ given by
\begin{eqnarray*}
\chi(b \cdot( f \odot \psi))_q(y) &=& \int_{\sigma \in \Theta^{\varrho(q)}} \int_{\tau \in \Xi^{\alpha(q)}} b(\tau)f((\tau^{-1} \cdot q) \cdot \sigma) (\varphi^{\#}_{\sigma^{-1}} \psi)_{\varrho(q)}  \\ 
&=& \int_{\tau \in \Xi^{\alpha(q)}} \int_{\sigma \in \Theta^{\varrho(q)}}  b(\tau)f(\tau^{-1} \cdot (q \cdot \sigma)) (\varphi^{\#}_{\sigma^{-1}} \psi)_{\varrho(q)} \\
&=& \int_{\tau \in \Xi^{\alpha(q)}} \int_{\sigma \in \Theta^{\varrho(q)}}  b(\tau)f(\tau^{-1} \cdot (q \cdot \sigma)) (\varphi^{\#}_{\sigma^{-1}} \psi)_{\varrho(\tau^{-1} \cdot q)} \\
&=& \int_{\tau \in \Xi^{\alpha(q)}} \int_{\sigma \in \Theta^{\varrho(q)}} b(\tau) \varphi_{\tau^{-1}}^{\#}(f(q \cdot \sigma) (\varphi^{\#}_{\sigma^{-1}} \psi))_{\varrho(q)} \\
&=& b \cdot \int_{\sigma \in \Theta^{\varrho(q)}}  f(q \cdot \sigma) (\varphi^{\#}_{\sigma^{-1}} \psi)_{\varrho(q)} \\
&=& b \cdot \chi(f \odot \psi)_q(y)  
\end{eqnarray*}
for all $y \in Y$ and $q \in \alpha^{-1}(y)$. The second equality holds because the left and the right actions of $K$ and $G$ on $Q$ are commutative. The third equality holds because the $K$ action on $Q$ preserves the fibres of $\varrho$. It follows that $\chi: \phi_{\#} \Gamma^{\infty}(\xi) \rightarrow \Gamma^{\infty}(\phi_{\#} \xi)$ is an isomorphism of $C^{\infty}(K \ltimes Y)$-modules.

The remaining task is to prove that the isomorphism $\chi$ extends to an isomorphism of Hilbert spaces. On $X$ we have the riemannian structure and the volume form $\nu$. We can use $\varrho$ to lift $\nu$ and the hermitian structure to $Q$. This makes $\varrho$ a riemannian covering map. Then 
\begin{eqnarray}\label{triangle}
\# K \int_X  (\psi, \psi') \nu = \pm \int_Q \varrho^*(\psi,\psi') \varrho^*(\nu) 
\end{eqnarray}
The $+1$ sign appears if $\varrho$ is orientation preserving, and $-1$ appears otherwise. To see that the equation holds, we give $Q$ a $G$-equivariant smooth triangulation which is fine enough so that $\varrho$ restricts to a diffeomorphism on each triangle, \cite{Ill78}. If $\triangle$ is any of these triangles in $Q$, then
\begin{eqnarray*}
\int_{\varrho(\triangle)} (\psi, \psi') \nu = \pm \int_{\triangle} \varrho^*(\psi, \psi') \varrho^* (\nu). 
\end{eqnarray*}
Now \eqref{triangle} holds because there are exactly $\# K$ triangles in $Q$ which map to $\varrho(\triangle)$ in $X$. Let $\phi_{\#} \nu$ denote the induced riemannian volume form on $Y$. By definition $\alpha^*(\phi_{\#} \nu) = \varrho^*(\nu)$. Therefore we have
\begin{eqnarray*}
\# G \int_Y  (\eta, \eta')_{\#} (\phi_{\#} \nu) &=& \pm \int_Q \alpha^*(\eta, \eta')_{\#} \alpha^*(\phi_{\#}\nu) \\
&=& \pm \int_Q \alpha^*(\eta, \eta')_{\#} \varrho^*(\nu) 
\end{eqnarray*}
for all $\eta$ and $\eta'$ in $\Gamma^{\infty}(\phi_{\#} \xi)$. Now the $+1$ sign appears if $\alpha$ is orientation preserving and $-1$ otherwise.

Suppose that $f \odot \psi$ and $f' \odot \psi'$ are in $C^{\infty}(Q) \odot_{\Theta} \Gamma^{\infty}(\xi)$, then 
\begin{eqnarray*}
& &\la f \odot \psi, f' \odot \psi' \ra  = \la (f',f)_{\Theta} \cdot \psi, \psi' \ra_X \\
&=& \int_X \int_{\Theta^{x}} \int_{\Xi^{\alpha(q)}}( \overline{f'(\tau^{-1} \cdot q)} f(\tau^{-1} \cdot q \cdot \sigma) \rho(\sigma) \psi_{\sigma^{-1} \cdot x}, \psi'_x) \nu \\
&=& \int_X \int_{\Theta^{x}} \int_{\Xi^{\alpha(q)}}(f(\tau^{-1} \cdot q \cdot \sigma) \rho(\sigma) \psi_{\sigma^{-1} \cdot x}, f'(\tau^{-1} \cdot q) \psi'_{x}) \nu \\
&=& \int_X \int_{\Theta^{(\sigma')^{-1} \cdot x}} \int_{\Theta^{x}} \int_{\Xi^{\alpha(q)}}(f(\tau^{-1} \cdot q \cdot \sigma' \cdot \sigma) \rho(\sigma' \cdot \sigma) \psi_{(\sigma' \sigma)^{-1} \cdot x}, f'(\tau^{-1} \cdot q \cdot \sigma') \rho(\sigma')\psi'_{(\sigma')^{-1} \cdot x}) \nu \\
&=& \int_X \int_{\Theta^{x}} \int_{\Theta^{x}} \int_{\Xi^{\alpha(q)}}(f(\tau^{-1} \cdot q \cdot \sigma'') \rho(\sigma'') \psi_{(\sigma'')^{-1} \cdot x}, f'(\tau^{-1} \cdot q \cdot \sigma') \rho(\sigma')\psi'_{(\sigma')^{-1} \cdot x}) \nu.
\end{eqnarray*}
where for each $x \in X$ the parameter $q$ is any of the points $\varrho^{-1}(x)$. The third equality holds by the $\Theta$-invariance of $(\cdot, \cdot)$ and the fourth by the left-invariance of the Haar measure. Moreover, on $Y$ we can fix any point $q \in \alpha^{-1}(y)$ for every $y \in Y$ and then write
\begin{eqnarray*}
( \chi(f \odot \psi), \chi(f' \odot \psi') )_{\#,y} = \int_{\Theta^{\varrho(q)}} \int_{\Theta^{\varrho(q)}}  ( f(q \cdot \sigma) \rho(\sigma) \psi_{\varrho(q \cdot \sigma )}, f'(q \cdot \sigma') \rho(\sigma') \psi'_{\varrho(q \cdot \sigma')} )
\end{eqnarray*}
In this notation we have
\begin{eqnarray*}
& &\la \chi(f \odot \psi), \chi(f' \odot \psi') \ra_Y \\ 
&=&\int_Y \int_{\Theta^{\varrho(q)}} \int_{\Theta^{\varrho(q)}} ( f(q \cdot \sigma) \rho(\sigma) \psi_{\varrho(q \cdot \sigma )}, f'(q \cdot \sigma') \rho(\sigma') \psi'_{\varrho(q \cdot \sigma')} ) (\phi_{\#} \nu) \\
&=&\int_Y \int_{\Theta^{\varrho(q)}} \int_{\Theta^{\varrho(q)}} \int_{\Xi^{y}} ( f(\tau^{-1} \cdot q \cdot \sigma) \rho(\sigma) \psi_{\varrho(\tau^{-1} \cdot q \cdot \sigma )}, f'(\tau^{-1} \cdot q \cdot \sigma') \rho(\sigma') \psi'_{\varrho(\tau^{-1} \cdot q \cdot \sigma')} ) (\phi_{\#} \nu) \\
&=&\int_Y \int_{\Theta^{\varrho(q)}} \int_{\Theta^{\varrho(q)}} \int_{\Xi^{y}} ( f(\tau^{-1} \cdot q \cdot \sigma) \rho(\sigma) \psi_{\varrho(q \cdot \sigma )}, f'(\tau^{-1} \cdot q \cdot \sigma') \rho(\sigma') \psi'_{\varrho(q \cdot \sigma')} ) (\phi_{\#} \nu) 
\end{eqnarray*}
The second equality holds by the $\Xi$-invariance of the induced pairing $(\cdot, \cdot)_{\#}$, and the third because the $\Xi$-action on $Q$ preserves the fibres of $\varrho$.  Together with the formulas in the previous paragraph this implies that
\begin{eqnarray*}
\# K \la f \odot \psi, f' \odot \psi' \ra = \pm \# G \la \chi(f \odot \psi), \chi(f' \odot \psi') \ra_Y
\end{eqnarray*}
for all $f \odot \psi$ and $f' \odot \psi'$ in $C^{\infty}(Q) \odot_{\Theta} \Gamma^{\infty}(\xi)$. Now the sign $+1$ appears if both $\varrho$ and $\alpha$ are either orientation preserving or orientation reversing, and $-1$ appears otherwise. Consequently, $\chi$ extends to an isomorphism of Hilbert spaces $\phi_{\#} L^2(\xi) \rightarrow  L^2(\phi_{\#} \xi)$. 

The elements of $C^{\infty}(K \ltimes Y)$ act as bounded operators on $C^{\infty}(Q) \odot_{\Theta} \Gamma^{\infty}(\xi)$ and on $\Gamma^{\infty}(\phi_{\#} \xi)$. Therefore these actions extend to make $\phi_{\#} L^2(\xi)$ and $L^2(\phi_{\#} \xi)$ representation spaces for $C^{\infty}(K \ltimes Y)$. Now $\chi$ is an isomorphism of representations. \5 $\square$

\section{Geometric Spectral Triple}

The goal of this section is to prove Theorem 1.

\3 \noindent \textbf{3.1.} Suppose that $G \ltimes X$ and $K \ltimes Y$ are Morita equivalent compact action orbifolds. Assume that $G \ltimes X$ is $G$-spin. This property is Morita invariant and it follows readily that $K \ltimes Y$ is $K$-spin. Let us fix a $G$-equivariant complex Dirac bundle on $X$. More precisely, we fix a $G$-invariant riemannian structure on $X$, a $G$-equivariant complex Dirac bundle on $X$ which is equipped with a unitary and Clifford compatible $G$-invariant connection and an inner product so that the unit vector fields in the Clifford bundle act as unitary transformations, \cite{LM89}. Let us then fix a Morita equivalence $\phi = (\alpha, Q, \varrho)$ between $G \ltimes X$ and $K \ltimes Y$. Then there is the spectral triple induced by $\phi$ which is denoted by, \cite{Har14b}
\begin{eqnarray*}
(C^{\infty}(K \ltimes Y), L^2(\phi_{\#} \xi), \phi_{\#}\eth, \phi_{\#}\omega). 
\end{eqnarray*}
The Hilbert space $L^2(\phi_{\#} \xi)$ is the completion of $\Gamma^{\infty}(\phi_{\#} \xi)$ in the $L^2$-inner product. The induced Dirac operator $\phi_{\#} \eth$ is defined globally over $Y$ by the relation
\begin{eqnarray*}
\eta_q(y) = \psi_{\varrho(q)} \5 \Rightarrow \5 (\phi_{\#}\eth \eta)_q(y) = (\eth\psi)_{\varrho(q)}.  
\end{eqnarray*}
There is also the following local description. If $y \in Y$, then for any $q \in \alpha^{-1}(y)$, there is a local section $\beta_q$ of $\alpha$ which maps $y$ to $q$ and the induced Dirac operator can be written by
\begin{eqnarray*}
(\phi_{\#} \eth)_q(y) = (\varrho \circ \beta_q)^* \circ \eth_{\varrho(q)} \circ ((\varrho \circ \beta_q)^{-1})^*
\end{eqnarray*}
This description is independent on the choice of the local section because the local sections have a unique germ at $y \in Y$. The induced Dirac operator satisfies the condition
\begin{eqnarray*}
\rho(\sigma)(\phi_{\#} \eth)_{q \cdot \sigma}(y)\rho(\sigma)^{-1} = (\phi_{\#} \eth)_q(y). 
\end{eqnarray*}
which makes it a well defined operator on $\Gamma^{\infty}(\phi_{\#} \xi)$. 

In 2.4 we introduced an isomorphism $\chi: \phi_{\#} L^2(\xi) \rightarrow L^2(\phi_{\#} \xi)$. We can use this isomorphism to define an unbounded operator $\widetilde{\phi_{\#} \eth}: \phi_{\#} L^2(\xi) \rightarrow \phi_{\#} L^2(\xi)$  by
\begin{eqnarray}\label{tilde-eth}
\widetilde{\phi_{\#} \eth} = \chi^{-1} \circ \phi_{\#} \eth \circ \chi.
\end{eqnarray}
The subspace $C^{\infty}(Q) \odot_{\Theta} \Gamma^{\infty}(\xi)$ is a dense domain for $\widetilde{\phi_{\#} \eth}$. Since $\varrho$ is a local diffeomorphism, there is a local section $\delta_{q}$ of $\varrho$ for each $q \in Q$ defined in a neighborhood of $\varrho(q)$ in $X$ so that $\delta_{q}(\varrho(q)) = q$. The germs of these local diffeomorphisms are unique. Let $\eth = \gamma(e_i^*) \nabla_{e_i}$ be the Dirac operator so that $e_i$ are linearly independent vector fields. In the following we give the operator \eqref{tilde-eth}  an explicit local realization. 

\3 \noindent \textbf{Proposition 5.} The operator $\widetilde{\phi_{\#} \eth}$ is given by
\begin{eqnarray}\label{prop5}
(\widetilde{\phi_{\#} \eth}) (f(q) \odot \psi_x) = \sum_{i = 1}^n \Big(\varrho^* (e_i(f \circ \delta_q))(q) \odot \gamma(e_i^*)\psi_x \Big) + f(q) \odot (\eth \psi)_x
\end{eqnarray}
in the dense subspace $C^{\infty}(Q) \odot_{\Theta} \Gamma^{\infty}(\xi)$ of $\phi_{\#} L^2(\xi)$. 

\3 \noindent Proof. For the sake of notation simplicity, we are using $\Theta = G \ltimes X$ again.

It is sufficient to show that the right hand side of \eqref{prop5} maps to $(\phi_{\#} \eth) \chi(f \odot \psi)$ under $\chi$. Let $y \in Y$ and $q \in \alpha^{-1}(y)$ be arbitrary and $x = \varrho(q)$. For each $q \cdot \sigma$ there is a local section $\delta_{q \cdot \sigma}$ of $\varrho$ defined in a neighborhood of $ \varrho(q \cdot \sigma) = \sigma^{-1} \cdot x$ so that $\delta_{q \cdot \sigma}(\sigma^{-1} \cdot x) = q \cdot \sigma$.  In a small enough neighborhood of $y \in Y$ we have
\begin{eqnarray*}
(\phi_{\#} \eth) \chi (f \odot \psi)_q (y) &=& \eth_{\varrho(q)} \Big( \int_{\Theta^{\varrho(q)}} f(q \cdot \sigma) (\varphi^{\#}_{\sigma^{-1}} \psi)_{\varrho(q)} \Big)  \\
&=& \eth_{x} \Big( \int_{\Theta^{x}} f \circ \delta_{q \cdot \sigma}(\sigma^{-1} \cdot x)  (\varphi^{\#}_{\sigma^{-1}} \psi)_{x} \Big) \\
&=& \eth_{x} \Big( \int_{\Theta^{x}} (\varphi^{\#}_{\sigma^{-1}} (f \circ \delta_{q }) \psi)_{x} \Big) \\
&=&   \int_{\Theta^{x}} (\varphi^{\#}_{\sigma^{-1}} \eth (f \circ \delta_{q }) \psi)_{x} \\
&=& \int_{\Theta^{x}} \Big((\varphi^{\#}_{\sigma^{-1}} (\eth(f \circ \delta_{q })) \psi  + (\varphi^{\#}_{\sigma^{-1}} (f \circ \delta_q) \eth \psi))_{x} \Big) \\
&=& \int_{\Theta^{x}}\sum_i^n \Big( (\varphi^{\#}_{\sigma^{-1}} e_i(f \circ \delta_q) \gamma(e_i^*)\psi)_x \Big) + \chi(f \odot \eth \psi)_q(y)  \\
&=&\int_{\Theta^{x}}\sum_i^n \Big( e_i(f \circ \delta_{q \cdot \sigma})(\sigma^{-1} \cdot x) (\varphi^{\#}_{\sigma^{-1}} \gamma(e_i^*)\psi)_x \Big) + \chi(f \odot \eth \psi)_q(y) \\
&=&\int_{\Theta^{\varrho(q)}}\sum_i^n \Big( \varrho^*(e_i(f \circ \delta_{q \cdot \sigma}))(q \cdot \sigma) (\varphi^{\#}_{\sigma^{-1}} \gamma(e_i^*)\psi)_{\varrho(q)} \Big) + \chi(f \odot \eth \psi)_q(y). 
\end{eqnarray*}
The fourth equality uses the $G$-invariance of $\eth$. Since the Haar measure is globally constant over $X$ there are no contributions from the vector fields acting on  the measure terms. \5 $\square$

\3 Suppose that the spectral triple on $C^{\infty}(G \ltimes X)$ is chirally graded. The chiral grading operator is a $G$-invariant operator and there is the induced chiral grading operator on $\phi_{\#} F_{\Sigma}$ which is defined by the relation
\begin{eqnarray*}
\eta_q(y) = \psi_{\varrho(q)} \5 \Rightarrow \5 (\phi_{\#} \omega)_q(y) \eta_q(y) = \omega_{\varrho(q)} \psi_{\varrho(q)}. 
\end{eqnarray*}
for all $y \in Y$ and $q \in \alpha^{-1}(y)$. Let $\widetilde{\phi_{\#} \omega} = \chi^{-1} \omega \chi$. The following proposition can be proved with the strategy of Proposition 5 but the details are more straightforward. 

\3 \noindent \textbf{Proposition 6.} The operator $\widetilde{\phi_{\#} \omega}$ is given by 
\begin{eqnarray*}
(\widetilde{\phi_{\#} \omega}) (f(q) \odot \psi_x) = f(q) \odot (\omega \psi)_x.  
\end{eqnarray*}
in the dense subspace $C^{\infty}(Q) \odot_{\Theta} \Gamma^{\infty}(\xi)$ of $\phi_{\#} L^2(\xi)$. 

\3 \noindent \textbf{3.2.} This subsection consists of a proof for Theorem 1 of 1.2. 

\3 Suppose that $G \ltimes X$ and $K \ltimes Y$ is the pair of compact action orbifolds introduced in 3.1. According to Proposition 3, $C^{\infty}(Q)$ is a pre-imprimitivity $C^{\infty}(K \ltimes Y)$ - $C^{\infty}(G \ltimes X)$-bimodule. So, with $E = C^{\infty}(Q)$, \textbf{M1} holds. We have the unitary equivalence of spectral triples
\begin{eqnarray*} 
(C^{\infty}(K \ltimes Y), L^2(\phi_{\#} F_{\Sigma}), \phi_{\#} \eth, \phi_{\#} \omega) \mapsto (C^{\infty}(K \ltimes Y), \phi_{\#} L^2(F_{\Sigma}), \widetilde{\phi_{\#} \eth}, \widetilde{\phi_{\#} \omega})
\end{eqnarray*}
which was described in Theorem 2 and 3.1. In particular, \textbf{M2} holds. It is a direct consequence of Proposition 5 in this manuscript and Theorem 2 in \cite{Har14b} that the operators 
\begin{eqnarray*}
T_{u} \circ \eth -  \widetilde{\phi_{\#} \eth} \circ T_{u} \5 \text{and} \5 T_u^* \circ \widetilde{\phi_{\#} \eth} - \eth \circ T^*_u 
\end{eqnarray*}
are bounded operators for all $u \in E$. Since $\eth$ and $\widetilde{\phi_{\#} \eth}$ are finitely summable operators as Dirac operators on a compact manifold, it follows by standard means that \textbf{M3} holds as well, see e.g \cite{GVF00}. The dimension condition \textbf{M4} holds because in the coordinates of $Y$, the induced spinor bundle $\phi_{\#} F_{\Sigma}$ and the induced Dirac operator $\phi_{\#} \eth$, are a actual spinor bundle and a Dirac operator on $Y$, \cite{Har14b} 4.3, and in this case the spectral dimension will be preserved because the dimension of the base is preserved. 

There is an isomorphism of vector spaces, $\phi_{\#}: \Gamma^{\infty}(F_{\Sigma})^G \rightarrow \Gamma^{\infty}(\phi_{\#} F_{\Sigma})^K$ (read $G$ and $K$-invariant subspaces), which can be defined as follows: for each $\psi \in \Gamma^{\infty}(F_{\Sigma})^G$, $(\phi_{\#} \psi)_q(y) = \psi_{\varrho(q)}$ where $q \in \alpha^{-1}(y)$. Similarly we have the isomorphism $\phi_{\#}: C^{\infty}(X)^G \rightarrow C^{\infty}(Y)^K$ which sends $f \in C^{\infty}(X)^G$ to $\phi_{\#}(f): y \mapsto f(\varrho(q))$ where $q \in \alpha^{-1}(y)$ is an arbitrary point. The choice is arbitrary by $G$-invariance of $f$. In \cite{Har14b} it was found that 
\begin{eqnarray}
\phi_{\#}(f) \phi_{\#}(\psi) = \phi_{\#}(f \psi)
\end{eqnarray}
for all $f \in C^{\infty}(X)^G$ and $\psi \in \Gamma^{\infty}(\xi)^G$, and 
\begin{eqnarray*}
\phi_{\#} \circ \eth \circ \phi_{\#}^{-1} = \phi_{\#} \eth, \5 \phi_{\#} \circ \omega \circ \phi_{\#}^{-1} = \phi_{\#} \omega. 
\end{eqnarray*}
As in the proof of Theorem 2, we can use the covering maps $\varrho$ and $\alpha$ to make $Q$ a riemannian covering space on $X$ and on $Y$. If $\psi, \psi' \in \Gamma^{\infty}(\xi)^G$ we have 
\begin{eqnarray*}
\# K \int_X (\psi, \psi') \nu  = \pm \int_Q \varrho^* (\psi, \psi') \varrho^*(\nu)
\end{eqnarray*}
where the sign $+1$ corresponds to the case of orientation preserving $\varrho$. Similarly we have the equality
\begin{eqnarray*}
\# G \int_Y (\phi_{\#} \psi, \phi_{\#} \psi')_{\#} \phi_{\#}\nu  &=& \pm \int_Q \alpha^*(\phi_{\#} \psi, \phi_{\#} \psi')_{\#} \alpha^*(\phi_{\#}\nu) \\
&=& \pm \int_Q \alpha^* \phi_{\#} (\psi, \psi') \alpha^*(\phi_{\#}\nu) \\
&=& \pm \int_Q \varrho^*(\psi, \psi') \varrho^*(\nu).
\end{eqnarray*}
where $+1$ corresponds to the case of orientation preserving $\alpha$. Consequently, 
\begin{eqnarray*}
\la \psi, \psi' \ra_X = \pm \frac{\# G}{\# K} \la \phi_{\#} \psi, \phi_{\#} \psi' \ra_Y. 
\end{eqnarray*}
holds for all $\psi, \psi' \in \Gamma^{\infty}(F_{\Sigma})^G$. The sign $+1$ appears in the case where both $\varrho$ and $\alpha$ are either orientation preserving or orientation reversing, and $-1$ appears otherwise. The linear map 
\begin{eqnarray*}
\sqrt{-1}^{\frac{1 \mp 1}{2}}\sqrt{\frac{\# G}{\# K}} \phi_{\#} : \Gamma^{\infty}(F_{\Sigma})^G \rightarrow \Gamma^{\infty}(\phi_{\#}F_{\Sigma})^K 
\end{eqnarray*}
extends to a unitary transformation $U_{\phi}: L^2(F_{\Sigma})^G \rightarrow L^2(\phi_{\#} F_{\Sigma})^K$ and $U_{\phi}$ determines a unitary equivalence of the invariant spectral triples
\begin{eqnarray*}
(C^{\infty}(X)^G, L^2(F_{\Sigma})^G, \eth, \omega) \rightarrow (C^{\infty}(Y)^K, L^2(\phi_{\#} F_{\Sigma})^K, \phi_{\#}\eth, \phi_{\#}\omega).
\end{eqnarray*}
So, \textbf{M5} holds. \5 $\square$

\3 In the references, \cite{Har14a} and \cite{Har14b}, the inner products in the invariant spectral triples are defined with respect to the integration theory over the groupoid orbit space $|\Theta|$ in the context of proper {\' e}tale Lie groupoids. However it is not hard to see that this inner product is equivalent to the $L^2$-inner on $X$ that was applied in this work when $\Theta = G \ltimes X$ is a compact action orbifold. More precisely, if $f \in C^{\infty}(X)^G$, then 
\begin{eqnarray*}
\# G \int_{| G \ltimes X| } f \nu = \int_X f \nu
\end{eqnarray*}
by similar triangulation statement that was used in the proof of Theorem 1 and Theorem 2. The effectiveness hypothesis for the $G$-action is essential here.

\section*{Appendix}

\noindent \textbf{A.1.} In the case of a smooth compact riemannian spin manifold $(X,g)$, there is a Dirac spectral triple associated with the metric $g$ and the geodesic distance is given in the spectral triple language by the formula
\begin{eqnarray*}
d_g(p,q) = \text{sup}\{|a(p)-a(q)|: a \in \mathcal{C}(X), ||[\eth,a]|| = ||\text{grad}(a)|| \leq 1 \}
\end{eqnarray*}
for all $p, q \in X$, \cite{Con94}. The supremum is taken over the space $\mathcal{C}(X)$ of all continuous functions whose gradient is defined almost everywhere as an essentially bounded measurable vector field. In this appendix, this theorem is specialized to compact action orbifolds with pointlike singularities. 

\3 \noindent \textbf{A.2.}  A piecewise smooth $G$-path from $x \in X$ to $x' \in X$ in the compact action orbifold $G \ltimes X$ is given by a subdivision $0 = t_0 < t_1 < \cdots < t_k = 1$ of $[0,1]$, and a sequence $\gamma = (c_1, \sigma_1, \ldots, \sigma_{k-1}, c_k)$ so that:
\begin{quote}
\textbf{1.} $c_i: [t_{i-1}, t_i] \rightarrow X$ is a piecewise smooth map for all $i \in \{1, \ldots, k\}$, and $c_1(0) = x$ and $c_k(1) = x'$, 

\textbf{2.} $\sigma_i \in G \ltimes X$ and $c_i(t_i) = \sigma_i \cdot c_{i+1}(t_i)$ for all $i \in \{1, \ldots, k-1\}$. 
\end{quote}
Suppose that a $G$-invariant riemannian structure has been fixed in $G \ltimes X$. The length of a $G$-path is defined by 
\begin{eqnarray}\label{lg}
l_g(\gamma) = \sum_{i = 1}^k l_g(c_i). 
\end{eqnarray}
The geodesic distance on the compact action orbifold is defined by
\begin{eqnarray*}
d_g^{X/G}(x,x') = \text{inf} \{ l_g(\gamma) \}
\end{eqnarray*}
where the infimum runs over all $G$-paths between $x$ and $x'$ in the orbifold, \cite{GH06}. Clearly this distance should be understood as the distance between the points $[x]$ and $[x']$ in the orbit space $X/G$. 

For the following theorem we suppose that the compact action orbifold $G \ltimes X$ is $G$-spin. Moreover, the $G$-invariant riemannian metric $g$, the $G$-equivariant spinor bundle $F_{\Sigma}$ and the $G$-invariant Dirac operator $\eth$ have been fixed. This data gives rise to the noncommutative orbifold determined by the Dirac spectral triple on $C^{\infty}(G \ltimes X)$. The associated invariant spectral triple $(C^{\infty}(X)^G, \hil^G, \eth)$ is studied in the following. 

\3 \noindent \textbf{Theorem 3.} Suppose that $G \ltimes X$ has either pointlike singularities or the $G$ action is free. The spectral distance formula of the invariant spectral triple computes the geodesic distances in the compact action orbifold $G \ltimes X$: 
\begin{eqnarray*}
d_g^{X/G}(x,x') = \text{sup}\{|a(x)-a(x')|: a \in \mathcal{C}(X)^G, ||[\eth,a]|| \leq 1 \},
\end{eqnarray*}
for all $x,x' \in X$. \3

\noindent Proof. Denote by $\Sigma_G$ the singular locus which consists of isolated points in $X$. We first study the case $x,x' \in X - \Sigma_G$. As an open subset of $X$, $X - \Sigma_G$ is a smooth manifold and $G$ acts on it freely. The metric $g$ restricts to define a $G$-invariant metric in $X - \Sigma_G$. The spaces $\mathcal{C}(X - \Sigma_G)^G$ and $\mathcal{C}(X)^G$ can be identified. Namely, if $a \in \mathcal{C}(X - \Sigma_G)$ then $a$ is a continuous function which has a bounded derivative almost everywhere and so $a$ extends uniquely to an element in $\mathcal{C}(X)$ since it is uniformly continuous. The inverse is the restriction of functions in $X$ to $X - \Sigma_G$. This identification passes to the $G$-invariant subspaces. Therefore we can write 
\begin{eqnarray*}
&&\text{sup}\{|a(x)-a(x')|: a \in \mathcal{C}(X)^G, ||[\eth,a]|| \leq 1 \} \\
&=& \sup \{|a(x) - a(x')| : a \in \mathcal{C}(X - \Sigma_G)^G, ||[\eth, a]|| \leq 1\}.
\end{eqnarray*}
We can view $X - \Sigma_G$ as a base space of an action groupoid under the free $G$-action. It is Morita equivalent to the unit groupoid on the smooth orbit space $(X - \Sigma_G)/G$. Let us denote by $\phi$ the Morita equivalence between $G \ltimes (X - \Sigma)$ and $1 \ltimes (X - \Sigma)/G$. We equip the orbit space with the induced metric $\phi_{\#}(g)$, the induced spinor bundle $\phi_{\#} F_{\Sigma}$ and the induced Dirac operator $\phi_{\#} \eth$. The spaces of continuous functions $C(X- \Sigma_G)^G$ and $C((X-\Sigma_G)/G)$ can be identified. Namely,  if $a \in C(X- \Sigma_G)^G$, then we can assign with $a$  the function $\phi_{\#}(a)$ on $(X - \Sigma_G)/G$ which sends $p$ to $a(x)$ if $x$ is any point so that $[x] = p$. The choice is arbitrary by the $G$-invariance of $a$. Now
\begin{eqnarray*}
||[\phi_{\#}\eth, \phi_{\#}(a)]|| &=& ||\phi_{\#}([\eth, a])|| \\
&=& \text{sup} \Big( \phi_{\#} ([\eth, a],[\eth, a])\Big)^{\frac{1}{2}} \\
&=&  \text{sup} \Big(([\eth, a],[\eth, a])\Big)^{\frac{1}{2}}
\end{eqnarray*}
for all $a \in C(X- \Sigma_G)^G$ which has the gradient defined almost everywhere. The first essential supremum is taken over $(X - \Sigma_G)/G$ and the second over $X - \Sigma_G$. It follows that $\phi_{\#}$ preserves the lengths of gradients. Now we have 
\begin{eqnarray}\label{dist}
&&  \sup \{|a(x) - a(x')| : a \in \mathcal{C}(X - \Sigma_G)^G, ||[\eth, a]|| \leq 1\} \\ \nonumber
&=& \sup \{|b([x]) - b([x'])| : b \in \mathcal{C}((X - \Sigma_G)/G), ||[\phi_{\#}\eth, b]|| \leq 1\}.
\end{eqnarray}
for all $x, x' \in X - \Sigma_G$. The right hand side of \eqref{dist} computes the geodesic distance on the manifold $(X - \Sigma_G)/G$ with respect to the induced metric $\phi_{\#}(g)$. This gives us us the equality 
\begin{eqnarray*}
\sup \{|b([x]) - b([x'])| : b \in \mathcal{C}((X - \Sigma_G)/G), ||[\phi_{\#}\eth, b]|| \leq 1\} = \inf \{ l_{\phi_{\#}(g)}(\gamma) \} 
\end{eqnarray*}
 where $l_{\phi_{\#}(g)}$ denotes the standard path length on the manifold $(X - \Sigma_G)/G$, and the infimum runs over all piecewise smooth paths in $(X - \Sigma)/G$ from $[x]$ to $[x']$.

If $\gamma$ is a piecewise smooth $G$-path in  $X - \Sigma_G$ then it determines a piecewise smooth path $[\gamma]$ in $(X - \Sigma_G)/G$ and the lengths satisfy $l_g(\gamma) = l_{\phi_{\#}(g)}([\gamma])$ where $l_g$ is defined as \eqref{lg} in the case of the orbifold $G \ltimes (X - \Sigma_G)$. Conversely, if $\gamma'$ is any piecewise smooth path in $(X - \Sigma_G)/G$ and if $\hat{\gamma}'$ is any $G$-path that is a lift of $\gamma'$ in $X - \Sigma_G$, then $l_g(\hat{\gamma}') = l_{\phi_{\#}(g)}(\gamma')$. These observations hold because with the riemannian metrics applied, $X - \Sigma_G \rightarrow (X - \Sigma_G)/G$ is a riemannian covering space. Consequently, 
\begin{eqnarray*}
\inf \{ l_{\phi_{\#}(g)}(\gamma') \} = \inf \{l_g(\gamma)\}
\end{eqnarray*}
where the first infimum runs over all piecewise smooth paths in $(X - \Sigma)/G$ from $[x]$ to $[x']$ and the second runs over all piecewise smooth $G$-paths in $X - \Sigma$ from $x$ to $x'$.

Suppose that dimension of $X$ is strictly greater than $1$. Then 
\begin{eqnarray*}
\inf \{l_g(\gamma)\} =  d_g^{X/G}(x,x')
\end{eqnarray*}
where the infimum runs over all piecewise smooth $G$-paths in $X - \Sigma$ from $x$ to $x'$. The right side denotes the orbifold geodesic distance on $G \ltimes X$. The equality holds because the geodesic distances are not affected if the set of isolated points $\Sigma_G$ is removed from $X$. In fact, the formula also holds trivially in the case where $X$ is one dimensional because all the compact action orbifolds with one dimensional base which are equipped with a spin structure (in particular they are orientable) have $\Sigma_G = \emptyset$. The string of equalities above gives a proof for the theorem in the special case of $x, x' \in X - \Sigma_G$.

Next assume that $x \in \Sigma_G$ and $x' \in X - \Sigma_G$. The geodesic distance function $d_g^{(X - \Sigma_G)/G}$ in the action groupoid $G \ltimes (X - \Sigma_G)$ defines a metric in its orbit space $(X - \Sigma_G)/G$ and therefore the function $X - \Sigma_G \rightarrow \mathbb{R}$ defined by 
\begin{eqnarray}\label{distance}
 z \mapsto d_g^{(X - \Sigma_G)/G}(z,x')
\end{eqnarray}
is uniformly continuous for all $z \in X - \Sigma_G$. Thus, \eqref{distance} extends uniquely to a continuous function in $X$ which sends $x$ to $d_g^{X/G}(x,x')$. On the other hand, according to the analysis above, the function \eqref{distance} is the same function as  
\begin{eqnarray*}
z \mapsto \text{sup}\{|a(z)-a(x')|: a \in \mathcal{C}(X - \Sigma_G)^G, ||[\eth,a]|| \leq 1 \}.
\end{eqnarray*}
in the domain $X - \Sigma_G$. This function extends to a continuous function over all $X$ which has to be equal to the extension of \eqref{distance} by uniqueness. So the theorem holds in this case. The analogous analysis shows that we can also take the end point of this path to be in the singular locus.  \5 $\square$

\end{document}